\documentclass[10pt]{amsart}
\usepackage{amssymb}

\usepackage{amscd}
\usepackage{amsthm}

\usepackage[all]{xy}
\usepackage[dvips]{graphicx}

\newtheorem{Thm}{Theorem}[section]
\newtheorem{Lem}[Thm]{Lemma}

\newtheorem{Prop}[Thm]{Proposition}

\newtheorem{Ques}[Thm]{Question}

\newenvironment{Proof}{\noindent{\it Proof.\ \
}}{\phantom{}\hfill$\square$}

\begin{document}

\title[Harder-Narasimhan filtrations and K-groups of an elliptic curve]
{Harder-Narasimhan filtrations and K-groups of an elliptic curve}

\author{Guodong ZHOU}
\address{Guodong ZHOU\newline
 LAMFA et CNRS UMR 6140 \newline Universit\'e de Picardie Jules Verne\newline
33, rue St Leu\newline 80039 Amiens\newline FRANCE}
\email{guodong.zhou@u-picardie.fr}


\begin{abstract}   Let $X$ be an elliptic curve over an algebraically closed field.
   We prove that some exact sub-categories of the category
of all  vector bundles over $X$,  defined using Harder-Narasimhan
filtrations,  have the same K-groups as the whole category.

\end{abstract}

\maketitle





\section{Introduction}
 Throughout this paper, $k$ denotes an algebraically closed field.
 Let $X$ be a smooth projective curve over $k$ and let $E$ be a
 vector bundle over $X$. We define the slope of $E$ as the quotient
 of its degree by its rank, i.e.
 $\mu(E)=deg(E)/rank(E)$.  A vector bundle $E$ is
 called semi-stable (resp. stable) if for any non-zero proper sub-bundle $E'$, we have
 $\mu(E')\leq \mu(E)$ (resp. $\mu(E')< \mu(E)$). The importance of
 the notion of semi-stability consists in the constructions of moduli
 spaces of vector bundles, see for example \cite{Mum}\cite{S}\cite{Sim}\cite{LP}\cite{HL}. For each vector bundle
 $E$, there exists a unique filtration, say \textit{Harder-Narasimhan
 filtration
 } (\cite[Proposition 1.3.9]{HN}),
 $$0=E_0\subsetneq  E_1\subsetneq \cdots \subsetneq E_{s-1}
 \subsetneq E_s=E$$
 such that the quotients $F_i=E_i/E_{i-1}$ are semi-stable for all
 $1\leq i\leq s$ and $$\mu(F_1)>\mu(F_2)>\cdots >\mu(F_s).$$
 We note  $\mu_{\text{max}}(E)=\mu(F_1)$ and
 $\mu_{\text{min}}(E)=\mu(F_s)$.

Let $\mathcal{P}(X)$ be the  exact category   of all vector bundles
over $X$.  Let $I\subset \mathbb{R}$ be a connected interval
(possibly of length zero). Following T. Bridgeland(\cite[Section
3]{B}), denote by $\mathcal{P}(I)$ the full sub-category of
$\mathcal{P}(X)$ consisting of all vector bundles $E$ such that
$\mu_{\text{max}}(E), \mu_{\text{min}}(E)\in I$. It is an
interesting fact
  that  the category $\mathcal{P}(I)$ is also exact with the exact
category structure induced from that of $\mathcal{P}(X)$ (see
Lemma~\ref{exact} below). We can therefore consider K-groups of
 $\mathcal{P}(I)$, as defined by D.Quillen for  an exact category
using his famous $Q$-construction (\cite{Q}). In this paper, we are
interested in the relations between K-groups of $\mathcal{P}(I)$ and
K-groups of $\mathcal{P}(X)$, i.e. those of $X$ in case that $X$ is
an elliptic curve.  More precisely, we prove the following theorem.
\begin{Thm}\label{main} Let   $X$ be an elliptic curve over $k$  an algebraically closed field
and let $I$ be a connected interval of strictly positive length.
Then the inclusion functor $\mathcal{P}(I) \hookrightarrow
\mathcal{P}(X)$ induces isomorphisms of K-groups
$K_i(\mathcal{P}(I)) \stackrel{\cong}{\rightarrow} K_i(X)$ for all
$i\geq 0$.
\end{Thm}

Vector bundles over an elliptic curve were classified by M.Atiyah in
\cite{A}. His classification is essential to the proof of the
preceding theorem. Roughly speaking, the idea    is to construct,
for an enough general vector bundle, a resolution of length one in
$\mathcal{P}(I)$ and then   the resolution theorem(\cite[Theorem
3.3]{Q}) applies.

The following question is natural.
\begin{Ques} Does the  statement in the preceding theorem hold if we replace $X$ by   any smooth
projective curve of genus $\geq 2$?
\end{Ques}

\textbf{Acknowledgement}:  I want to express my   gratitude to Mr.
Bruno Kahn (Directeur de Recherche du CNRS)  for some useful
conversations.


 \section{Proof of the main theorem}

Firstly we prove   the following fact mentioned in the introduction.

\begin{Lem} \label{exact}The category $\mathcal{P}(I)$ is an exact category
whose exact sequences are given by short exact sequences in
$\mathcal{P}(X)$ with their terms   in $\mathcal{P}(I)$.

\end{Lem}

\begin{Proof} One needs to show  that $\mathcal{P}(I)$ is closed under
extensions.
 Take a short exact sequence $$0\rightarrow E' \rightarrow E\rightarrow
E''\rightarrow 0$$ with $E', E''\in \mathcal{P}(I)$. Let
$$0=E_0\subsetneq  E_1\subsetneq \cdots \subsetneq E_{s-1}
 \subsetneq E_s=E$$ be the Harder-Narasimhan filtration of $E$.
 We then have  the exact sequence
$$0\rightarrow E'\cap E_{s-1} \rightarrow E'
 \rightarrow F \rightarrow 0$$ with $F$ a sub-bundle of
 $F_{s-1}=E/E_{s-1}$.  We obtain that
 $$\mu_{\text{min}}(E)=\mu(F_{s-1}) \geq \mu(F)\geq
 \mu_{\text{min}}(E').$$ We also have the exact sequence
$$0\rightarrow E'\cap E_{1} \rightarrow E_1
 \rightarrow G \rightarrow 0$$ with $G$ a sub-bundle of
 $E''$.  We get that
 $\mu(E'\cap E_1)\leq \mu_{\text{max}}(E')$ and $ \mu(G)\leq \mu(E'')
 \leq
 \mu_{\text{max}}(E'')$ and as $\mu_{\text{max}}(E)=\mu(E_1)$ is the
 barycenter  of $\mu(E'\cap E_1)$ and $\mu(G)$ with positive coefficients, $\mu_{max}(E)\leq \mu_{max}(E'), \mu_{max}(E'')$.
 This prove that $E\in \mathcal{P}(I)$.

\end{Proof}

\bigskip

Next we recall   some   known facts about vector bundles over an
elliptic curve  $X$.

\begin{Lem}\cite[Chapter 8, Section 8.7, Exercise 2.2]{LP}
  Each vector bundle over $X$ is a direct sum of
  indecomposable bundles. In particular, every indecomposable vector
  bundle is semi-stable.
  \end{Lem}

\begin{Thm}\label{produit} Let $E$ and $F$ be two semi-stable vector bundles over
  an elliptic curve. Then $E  \otimes F$ is still semi-stable.
\end{Thm}

In fact, in case of   characteristic zero, the tensor product of two
semi-stable vector bundles is semi-stable over a smooth projective
curve of arbitrary genus. This was first proved by M.S.Narasimhan
and C.S.Seshadri using analytic method (\cite{NS})  and then by Y.
Miyaoka using algebraic method (\cite[Corollary 3.7]{M}). The case
of positive characteristic uses the notion of strong semi-stability.
A vector bundle is called strongly semi-stable if all its Frobenius
pullbacks are semi-stable. T. Oda proved in \cite[Theorem 2.16]{O}
(see also \cite[Corollary $3^p$]{SB}) that a semi-stable vector
bundle over an elliptic curve is strongly semi-stable. Then the
preceding theorem follows form the facts (\cite[Section 5]{M}) that
the tensor product of two strongly semi-stable vector bundles is
still strongly semi-stable and that strong semi-stability implies
semi-stability.

Let $\mathcal{E}(r, d)$ with $r\geq 1$ and $d\in \mathbb{Z}$ be the
set of isomorphism classes of indecomposable vector bundles of rank
$r$ and of degree $d$. When $r$ and $d$ are coprime, M. Atiyah
introduced a distinguished vector bundle $E_{r,d}\in
\mathcal{E}(r,d)$ (Atiyah noted it by $E_A(r, d)$) with the property
$E_{r, d}^*\cong E_{r, -d}$ (\cite[Corollary of Theorem 7]{A}).

\bigskip

Let us    construct the resolutions of length one for an enough
general vector bundle.  The starting point is the following lemma.

\begin{Lem}\label{resl}  Let $E\in \mathcal{E}(r,d)$ with $r\geq 1$ and $d>0$. Then there
exists a vector bundle $E'\in \mathcal{E}(r+d,d)$, unique up to
isomorphisms, given by the extension
$$0\rightarrow H^0(E)\otimes \mathcal{O}_X \rightarrow E'\rightarrow
E\rightarrow 0$$ Moreover, $H^0(E)\cong H^0(E')$ and the map
$H^0(E')\otimes \mathcal{O}_X \cong H^0(E)\otimes \mathcal{O}_X
\rightarrow E'$ is the evaluation map.

\end{Lem}

\begin{Proof} The existence of $E'$ follows from
\cite[Lemma 16]{A} and  other statements are easy consequences of
\cite[Lemma 15]{A}.
\end{Proof}

\begin{Prop}\label{resolution} Let $E\in \mathcal{E}(r, d)$ with $r\geq 1$ and $d>0$
and let $\epsilon >0$. There exists a short exact sequence
$$0\rightarrow E_1\rightarrow E_0 \rightarrow E\rightarrow 0$$ where  $E_1$ is
semi-stable of zero slope and where $E_0$ is semi-stable of slope  $
 \mu(E_0)\in (0, \epsilon) $.

\end{Prop}

\begin{Proof}
The preceding lemma gives an exact sequence
$$0\rightarrow H^0(F_1)\otimes \mathcal{O}_X \stackrel{ ev}{\rightarrow}
F_1\stackrel{f_1}{\rightarrow} E\rightarrow 0$$ with $F_1\in
\mathcal{E}(r+d, r)$ and where $ ev$ is the evaluation map. We again
 apply
 Lemma~\ref{resl} to $F_1$ and we obtain
$$0\rightarrow H^0(F_2)\otimes \mathcal{O}_X \stackrel{ ev }{\rightarrow}
F_2\stackrel{f_2}{\rightarrow} F_1\rightarrow 0$$ with $F_2\in
\mathcal{E}(r+2d,d)$. These two exact sequences yield
$$0\rightarrow \text{Ker}(f_1\circ f_2) \rightarrow F_2
 \stackrel{f_1\circ f_2}{\rightarrow} E\rightarrow 0$$ and
 $$0\rightarrow H^0(F_2)\otimes \mathcal{O}_X  \rightarrow \text{Ker}(f_1\circ f_2)
 \rightarrow  H^0(F_1)\otimes \mathcal{O}_X\rightarrow 0 \ \  (*)$$
   Lemma~\ref{exact} implies that $\text{Ker}(f_1\circ f_2)$ is semi-stable of
 zero slope.


 If we iterate this process for $n$ times with $n$ enough great such that $d/(r+nd)<\epsilon$, we
 get
 $$0\rightarrow \text{Ker}(f_{1}\circ \cdots \circ f_n) \rightarrow
 F_n
 \stackrel{f_{1}\circ \cdots \circ f_n}{\rightarrow} E\rightarrow 0 \ \
 (**).$$
 As above, it is easy to show that  $\text{Ker}(f_n\circ \cdots \circ f_0)$ is
 semi-stable of zero slope  and that (**) is the desired resolution.

 \end{Proof}

\bigskip

Now we give   the proof of the main theorem.

\bigskip

\begin{Proof}(of Theorem~\ref{main})

We can suppose that $I=(a, b)$ with $-\infty <a<b<+\infty$. For any
real number $\lambda$, we note $I+\lambda=(a+\lambda, b+\lambda)$.
Set $J=(a, +\infty)$.

Step I: \textit{We   show that the inclusion functor
$\mathcal{P}(I)\hookrightarrow \mathcal{P}(J)$ induces  isomorphisms
of K-groups.} Take two integers $r\geq 1$ and $d$ such that
$-\frac{d}{r}=-\mu \in I$,
 $(r,d)=1$ and $(r,p)=1$ if $\text{char} k=p>0$. By
 Theorem~\ref{produit}, the tensor product by $E_{r,d}$ is an exact
 functor from $\mathcal{P}(I)$ to $\mathcal{P}(I+\mu)$.
 Note that $0\in I+\mu$.
Let $E\in \mathcal{P}(J)$. Then $E\otimes E_{r,d} \in
 \mathcal{P}((a+\mu,
 +\infty))$. Suppose that $E\otimes E_{r,d}=\oplus F_i $ with all $F_i$
 indecomposable. If $F_i \in \mathcal{P}(I+\mu)$, then we
 take the resolution $$0\rightarrow 0\rightarrow
 F_i\stackrel{Id}{\rightarrow} F_i\rightarrow 0$$ and if
 $F_i \not \in \mathcal{P}(I+\mu)$, we take the resolution
  given by  Proposition~\ref{resolution} with $\epsilon=b+\mu$. The sum of these
  resolutions of all $F_i$ is a
  resolution of $E$ of the form
$$0\rightarrow E_1\rightarrow E_0 \rightarrow
 E\otimes E_{r,d}\rightarrow 0$$
 where $E_1$ is semi-stable of zero slope and where $E_0$ is in
 $\mathcal{P}(I+\mu)$.
Now the tensor product of  the   resolution above  by $E_{r, -d}$
gives
 $$0\rightarrow E_1\otimes E_{r, -d} \rightarrow E_0 \otimes E_{r, -d}
 \stackrel{f}{\rightarrow}
 E\otimes E_{r,d}\otimes E_{r, -d}\rightarrow 0$$
 By \cite[Corollary 2.7]{O},  $E_{r,d}\otimes E_{r, -d}\cong
\mathcal{E}nd(E_{r,d})=\mathcal{O}_X\oplus  G$. We write $g$ the
projection from $E\otimes E_{r,d}\otimes E_{r, -d}$ to $E\otimes G$.
We have an exact sequence
$$0\rightarrow E_1\otimes E_{r, -d} \rightarrow \text{Ker}(g\circ f)
  \rightarrow
 E \rightarrow 0$$
 Obviously  $E_1\otimes E_{r,-d}$ is semi-stable of slope $-\mu$.
 The inequality $\mu_{max}(\text{Ker}(g\circ f))\leq \mu_{max}(E_0\otimes E_{r, -d})$
 together with Lemma~\ref{exact}
 implies that
 $\text{Ker}(g\circ f) \in \mathcal{P}(I)$.
 The resolution theorem applies and we obtain that the inclusion functor
 $\mathcal{P}(I) \hookrightarrow \mathcal{P}(J)$ induces
 isomorphisms of K-groups.

Step II: \textit{We  show that the inclusion functor
$\mathcal{P}(J)\hookrightarrow \mathcal{P}(X)$ induces  isomorphisms
of K-groups.}
 By a theorem of Serre (\cite[Chapter 2,
Theorem 5.17]{H}), for each $E\in \mathcal{P}(X)$, we have an exact
sequence $$0\rightarrow E  \rightarrow \mathcal{O}_X(n)^m
  \rightarrow
 F \rightarrow 0$$ with $n, m>>0$ and then $\mathcal{O}_X(n)^m, F\in \mathcal{P}(J)$.
  Let us  consider the functor $\mathcal{P}(J)^{op} \hookrightarrow
 \mathcal{P}(X)^{op}$ where $op$ means the opposite category.
 Notice that  $Q\mathcal{C}^{op}\cong Q\mathcal{C}$(\cite[Page 94]{Q}) where $Q$ is the $Q$-construction and then
 $K_i(\mathcal{C}^{op})\cong K_i(\mathcal{C})$ for all $i\geq 0$,
 we can deduce from the resolution
 theorem that the inclusion functor
 $\mathcal{P}(J ) \hookrightarrow \mathcal{P}(X)$ induces
 isomorphisms of K-groups.

This finishes the proof.

\end{Proof}

\end{document}